\documentclass[sts]{imsart}

\usepackage{amssymb,amsbsy,amsmath,amscd,amsthm,amsfonts,MnSymbol}
\usepackage{cite}
\usepackage[utf8]{inputenc}
\usepackage{enumerate}
\usepackage{hyperref}
\usepackage{dsfont}
\usepackage{graphicx}
\usepackage{enumitem}
\usepackage{wrapfig}

\usepackage{algorithm}
\usepackage{algorithmic}

\usepackage{enumerate}

\newtheorem{definition}{Definition}
\newtheorem{theorem}{Theorem}

\newtheorem{example}{Example}
\newtheorem{lemma}{Lemma}

\newtheorem{assumption}{Assumption}

\newcommand{\R}{\mathbb R}
\newcommand{\E}{\mathbb E}

\newcommand{\dist}{\mathrm{d_H}}

\newcommand{\LS}{\mathrm{LS}}
\newcommand{\GS}{\mathrm{GS}}
\newcommand{\Sm}{\mathrm{S}}

\newcommand{\Var}{\textrm{Var}}

\newcommand{\DS}{\displaystyle}

\newcommand{\EE}{\mathbb{E}}

\newcommand{\PP}{\mathbb{P}}

\newcommand{\RR}{\mathbb{R}}

\begin{document}

\begin{frontmatter}

\title{Differentially private sub-Gaussian location estimators}
\runtitle{DP location estimators}



\author{Marco Avella-Medina \thanks{Columbia University, Department of Statistics, marco.avella@columbia.edu} \and Victor-Emmanuel Brunel \thanks{ENSAE ParisTech, Department of Statistics, victor.emmanuel.brunel@ensae.fr}}

\runauthor{M. Avella-Medina and V.-E. Brunel}

\setattribute{abstractname}{skip} {{\bf Abstract:} }

\begin{abstract}
 We tackle the problem of estimating a location parameter with differential privacy guarantees and sub-Gaussian deviations. Recent work in statistics has focused on the study of estimators that achieve sub-Gaussian type deviations even for heavy tailed data. We revisit some of these estimators through the lens of differential privacy and show that a naive application of the Laplace mechanism can lead to sub-optimal results. We design two private algorithms for estimating the median that lead to estimators with sub-Gaussian type errors. Unlike most existing differentially private median estimators, both algorithms are well defined for unbounded random variables that are not even required to have finite moments. We then turn to the problem of sub-Gaussian mean estimation and show that under heavy tails natural differentially private alternatives lead to strictly worse deviations than their non-private sub-Gaussian counterparts. This is in sharp contrast with recent results that show that from an asymptotic perspective the cost of differential privacy is negligible.

\end{abstract}

\begin{keyword}
\kwd{Differential Privacy}
\kwd{Location estimators}
\kwd{Sub-Gaussian estimators}
\kwd{Median}
\end{keyword}

\end{frontmatter}

\section{Introduction}

Differential privacy has emerged as a rigorous mathematical approach to privacy that has been extensively studied in the theoretical computer science and machine learning literature following the path breaking work of \cite{dworketal2006}.  In this framework one assumes that there is a trusted curator that holds some data containing some possibly sensitive records of $n$ individuals. The goal of privacy is to simultaneously protect every individual record while releasing global characteristics of the database\cite{dworkandroth2014}.

Even though the machine learning community has been very prolific in developing differentially private algorithms for complex settings such as multi-armed bandit problems \cite{mishraandthakurta2015, tossouandimitrakakis2016,shariffandsheffet2018},  high-dimensional regression \cite{kiferetal2012,talwaretal2015} and  deep learning \cite{abadietal2016,lecuyeretal2018},  some  basic statistical questions are only starting to be understood. For example, the first statistical minimax rates of convergence under differential privacy were recently established in \cite{duchietal2018,caietal2019}. Some earlier work framing differential privacy in traditional statistics terms include  \cite{wassermanandzhou2010, lei2011,smith2011,chaudhuriandhsu2012,karwaandslavkovic2016}.  Recent work  has also sought to develop differential privacy tools for statistical inference and hypothesis testing \cite{gaboardi2016,sheffet2017,avella2019,barrientosetal2019}

In this paper we revisit the simple statistical problem of location parameter estimation and study the non-asymptotic deviations of differentially private location parameter estimators. More specifically, we consider the problem of constructing median and mean estimators that achieve sub-Gaussian deviations under heavy tails.

\subsection{Motivation}
It is well known that given a random iid sample $X_1,\dots,X_n$ of sub-Gaussian random variables with $\EE[X_1]=\mu$ and $\Var[X_1]=\sigma^2$,  the empirical mean $\overline{X}_n=\frac{1}{n}\sum_{i=1}^nX_i$ satisfies with probability at least $1-\alpha$
\begin{equation}\label{mean}
 |\overline{X}_n-\mu|\leq \sqrt{\frac{2\sigma^2\log(2/\alpha)}{n}}.
\end{equation}
The accuracy of the empirical mean estimator expressed in the above deviation inequality is a direct consequence of the sub-Gaussian assumption. In fact the dependence on $2/\alpha$ in the error worsens significantly when the distribution does not have a moment generating function. In particular, when $X_1$ is only assumed to have two finite moments, one cannot get an error whose order is smaller than $\sqrt{2\sigma^2/(n\alpha)}$ as shown in \cite{catoni2012}.  
For the empirical median on the other hand, one does not even need to assume any finite moments in order to establish similar sub-Gaussian deviations.

In light of \cite{caietal2019}, one may naturally wonder how differential privacy will affect the deviation bounds discussed above. The statistical minimax rates established in \cite{caietal2019} show that the rates of convergence of differentially private mean estimators are described by two terms. The first one correspond to the usual parametric $1/\sqrt{n}$ convergence, while the second one is driven by the differential privacy tuning parameters and  and is of the order $1/n$. Consequently for large $n$ differential privacy does not come at the expense of slower statistical convergence rates. Even though it seems intuitive that a similar phenomenon from a non-asymptotic deviations perspective deviations there are several technical obstacles that need to be addressed before one can attempt to establish such results.

A notorious technical difficulty that renders the study of non-asymptotic deviations challenging for differentially private mean and median estimators is the vast majority of existing algorithms requires the input data to be bounded; see for example \cite{dworketal2006,nissimetal2007,lei2011,smith2011,bassilyetal2014}. This is clearly unsatisfactory from a theoretical and practical perspective as it rules out common distribution used in statistical modelling such as the normal, gamma and t-distributions just to name a few. This is particularly disturbing for median estimators since the usual non-private empirical median does not even need the existence of finite moments in order to exhibit sub-Gaussian deviations.

\subsection{Our contributions}

\begin{itemize}
\item We derive the first differentially private  median estimators with sub-Gaussian  type errors under minimal  conditions. This should arguably be one of the simplest differentially private estimators for which one can establish sub-Gaussian deviations without assuming bounded variables or the existence of a moment generating function. We attain this objective for two new differentially private median estimators that build on two different carefully calibrated variants of the popular Laplace mechanism. Indeed, the naive Laplace mechanism is not directly applicable to the median when the support of the data is unbounded. The first algorithm adapts the idea of smooth sensitivity calibration introduced in \cite{nissimetal2007} for our median estimation problem with unbounded data without any finite moment assumptions. This first construction requires a truncation step of the data that can be avoided by our second algorithm. The latter revisits the propose-test-release paradigm introduced in \cite{dworkandlei2009} and leads to the desired sub-Gaussian deviations while avoiding any truncation. This is achieved by carefully controlling for the occurrence of unfavorable data configurations that are observed with negligible probability. 

We would like to highlight that two other mechanisms, essentially tailored for the estimation of location parameters,  have also been studied: exponential mechanism \cite{chaudhuriandhsu2012} and perturbation of the loss function for $M$-estimators \cite{chaudhurietal2011}. However, these mechanisms seem not to be applicable to either the mean (unless the data are bounded), or the median. Indeed, in the case of the mean, the loss function is not Lipschitz and in the case of the median, the loss function is not smooth, which are two requirements for the validity of these mechanisms. Therefore, it seems that even the seemingly simple problem of estimating the median of unbounded random variables in a differentially private fashion is more involved than one would expect. 

\item


We explore the possibility of constructing differentially private mean estimators with sub-Gaussian type leading errors terms when the underlying distributions are only assumed to have two finite moments.  This question is motivated by the significant attention given to the study of sub-Gaussian mean estimators over the last years  \cite{catoni2012, bubecketal2013, devroyeetal2016} and the numerous successful extensions of these methods to more complex models including \cite{hsuandsabato2016,minsker2018,avellaetal2018,lugosiandmendelson2019,lecueandlerasle2019}.  Unlike our private median estimators, natural differential private counterparts of known sub-Gaussian mean estimators fail to yield the desired deviations. This  suggests that differentially private mean estimation might come at the expense of worse non-asymptotic high probability deviations.

\end{itemize}

\section{Preliminaries}

\subsection{Definitions}

For $x=(x_1,\ldots,x_n)\in\R^n$, we denote by $x_{(1)}, \ldots, x_{(n)}$ the reordered coordinates of $x$ in nondecreasing order, i.e. $\DS \min_{1\leq i\leq n} x_i=x_{(1)}\leq \ldots\leq x_{(n)}=\max_{1\leq i\leq n} x_i$. We let $\ell=\lfloor n/2\rfloor$ and $\hat m(x)=x_{(l)}$ be the empirical (left) median of $x$.

For $T>0$ and $u\in\R$, let the truncation operator $f_T$ be defined as $f_T(u)=u$ if $|u|\leq T$, $f_T(u)=\mathrm{sign}(u)T$ otherwise. For $x\in\R^n$, let $\hat m_T(x)=\hat m(y)$, where $y=(y_1,\ldots,y_n)$ is defined as the truncated version of $x$ at level $T$, i.e., $y_i=f_T(x_i), i=1,\ldots,n$.

For any two vectors $x,x'\in\R^n$, we define their Hamming distance $\dist(x,x')$ as the number of coordinates that differ in $x$ and $x'$: $\DS \dist(x,x')=\#\{i=1,\ldots,n:x_i\neq x_i'\}$, where $\#$ stands for cardinality. Before formally defining $(\varepsilon,\delta)$-differential privacy, we first review some useful notions of sensitivity to the data that will serve as building blocks in the construction of our algorithms.


\begin{definition}
	Let $h:\R^n\to\R$ be a given function. 
	\begin{enumerate}
		\item The local sensitivity of $h$ maps any data point $x\in\R^n$ to the (possibly infinite) number $$\LS_h(x)=\sup\{|h(x')-h(x)|: x'\in\R^n, \dist(x,x')\leq 1\}.$$
		 
		\item The global sensitivity of $h$ is the (possibly infinite) number $$\GS_h=\sup_{x\in\R^n}\LS_h(x).$$ 
		\item For all $\beta>0$, the $\beta$-smooth sensitivity of $h$ is the mapping 
	$$\DS \Sm_h^{(\beta)}(x)=\sup_{x'\in\R^n}\left(e^{-\beta \dist(x,x')}LS_{h}(x')\right), \quad x\in\R^n.$$
	\end{enumerate}
\end{definition}

\begin{example}
\begin{enumerate}
	\item It is easy to see that for all $x\in\R^n$, the local sensitivity of the empirical median is
	$$\LS_{\hat m}(x)=\max\left(x_{(\ell+1)}-x_{(\ell)},x_{(\ell)}-x_{(\ell-1)}\right).$$
	Moreover, for all $\beta>0$ and all $x\in\R^n$, $\DS \Sm_{\hat m}^{(\beta)}(x)=\infty$, and $\GS_{\hat m}=\infty$.
	\item Let $\DS \hat \mu(x)=n^{-1}\sum_{i=1}^n x_i, x\in\R^n$ be the empirical mean function. Then, all the above quantities are infinite.
	\item Let $\DS \hat\mu_T(x)=n^{-1}\sum_{i=1}^n f_T(x_i)$ be the empirical mean of the truncated entries of $x\in\R^n$. Then, for all $x\in\R^n$, $\DS \LS_{\hat\mu_T}(x)=n^{-1}\max\left(T-f_T(x_{(1)}),f_T(x_{(n)})+T\right)$. Thus, $\GS_{\hat\mu_T}=2T$ and $\Sm_{\hat\mu_T}(x)\leq 2T$, for all $x\in\R^d$. 
\end{enumerate}	
\end{example}

Let us now compute the smooth sensitivity of the empirical median of truncated numbers.

\begin{lemma} \label{lemma:comput-LS}
	Let $T>0$ and $x\in\R^n$. Let $y=(y_1,\ldots,y_n)$ with $y_i=f_T(x_i), i=1,\ldots,n$. For all $k<0$ (resp. $k>n$), set $y_{(k)}=-T$ (resp. $y_{(k)}=T$). Then, for all $\beta>0$, $\DS \Sm_{\hat m_T}^{(\beta)}(x)=\max_{k\geq 0}e^{-\beta k}\max_{t=0,\ldots,k+1} \left(y_{(\ell+t)}-y_{(\ell+t-k-1)}\right)$.
\end{lemma}

This lemma is a direct consequence of \cite{nissimetal2007}. The same authors also showed that $\DS \Sm_{\hat m_T}^{(\beta)}(x)$ can be computed in $O(n\log n)$ operations .

In what follows, we refer as \textit{random function} to any function $\tilde h:\R^n\to\R$ such that for all $x\in\R^n$, $\tilde h(x)$ is a Borelian random variable. In this paper, we will use the hat sign to denote non-randomized estimators, and the tilde sign to denote their randomized version (e.g., $\hat h$ vs. $\tilde h$).

\begin{definition}
	Let $\varepsilon,\delta>0$. A random function $\tilde h$ is called $(\varepsilon,\delta)$-differentially private if and only if for each pair $x,x'\in\R^n$ with $\dist(x,x')\leq 1$ and for all Borel sets $B\subseteq \R$, 
$$\PP[\tilde h(x)\in B]\leq e^\varepsilon \PP[\tilde h(x')\in B]+\delta.$$
\end{definition}

It is important to note that for all $\beta\geq \beta'>0$, all functions $h:\R^n\to\R$ and all $x\in\R^n$, it holds that
\begin{equation} \label{eq:order-sens}
	\LS_h(x)\leq \Sm_h^{(\beta)}(x)\leq \Sm_h^{(\beta')}(x)\leq \GS_h(x)\leq \infty.
\end{equation}
In particular, using the global sensitivity is more restrictive than using the smooth sensitivities, which itself is more restrictive than using the local sensitivity.

\subsection{Background}

The Laplace mechanism is one of the basic tools used in the differential privacy literature in order to construct private algorithms. The basic idea is to make deterministic functions private by adding  random noise calibrated using their sensitivity to the data. Let us review some well know results that we will use in the construction of our differentially private estimators.
Recall that the Laplace distribution with parameter $\lambda>0$ is the continuous probability distribution with density $(\lambda/2)e^{-\lambda|u|}, u\in\R$. We denote this distribution by $\mathrm{Lap}(\lambda)$. 

The following theorem, due to \cite{dworketal2006}, gives a very simple way to make a function $h$ differentially private. However, it requires the very strong assumption that $h$ has a finite global sensitivity.

\begin{theorem} \label{thm:DP-GS}
	Let $h:\R^n\to\R$ be a function with finite global sensitivity. Let $Z$ be a Laplace random variable with parameter $1$. For all $\varepsilon>0$, the random function $\DS \tilde h(x)=h(x)+\frac{Z}{\varepsilon}\GS_h, x\in\R^n$, is $(\varepsilon,0)$-differentially private.
\end{theorem}

The following result is due to \cite{nissimetal2007}, and allows for less restrictive functions $h$.

\begin{theorem} \label{thm:DP-SS}
	Let $h:\R^n\to\R$ and assume that for all $k=1,\ldots,n$ and all $x\in\R^n$, $\LS_h^{(k)}(x)<\infty$. Let $Z$ be a Laplace random variable with parameter $1$. Let $\varepsilon,\delta>0$ and set $\DS \beta=\frac{\varepsilon}{2\log(1/\delta)}$. Then, the random function $\DS \tilde h(x)=h(x)+\frac{2Z}{\varepsilon}\Sm_h^{(\beta)}(x), x\in\R^n$, is $(\varepsilon,\delta)$-differentially private.
\end{theorem}

Note that these two versions of the Laplace mechanisms, using either the global, or the smooth sensitivities, cannot be used directly for the empirical mean or the empirical median when the data are unbounded, since the two sensitivities are infinite. 
In the next section, we first adapt the Laplace mechanism with the smooth sensitivity, by truncating the data in order to ensure a finite smooth sensitivity. Note that the local sensitivity, which is the least restrictive in the sense of \eqref{eq:order-sens}, cannot be used directly in a Laplace mechanism, as explained in \cite{nissimetal2007}. In the second part of next section, we define an estimator that uses some sort of local sensitivity, which is less restrictive than the smooth sensitivity and thanks to which we can avoid truncating the data.

\section{Median estimation}
\label{medians}

In this section we describe two differentially private median estimators. As explained above, the first one uses the Laplace mechanism but forces us to truncate the data, and to assume that the true (and unknown) median of the data is bounded. In the second algorithm, we no longer truncate the data, by using a more subtle approach than the Laplace mechanism. Our differentially private algorithms seem to be the first ones to achieve optimal statistical rates of convergence under no moment assumptions on the data.

\subsection{Private median via smooth sensitivity calibration} \label{sec:MedianSS}

We will only require the following distributional assumption in the derivation of our deviation inequalities for our  private median estimators.

\begin{assumption} \label{Ass:1}
	The distribution of $X_1$ has a density $f$ with respect to the Lebesgue measure and it has a unique median $m$. Moreover, there exist positive constants $r,L$ such that $f(u)\geq L$, for all $u\in [m-r,m+r]$.
\end{assumption}

In particular, under this assumption, the cdf $F$ of $X_1$ satisfies the following: 
\begin{equation} \label{comment-ass}
	|F(u)-F(v)|\geq L|u-v|, \forall u,v\in [m-r,m+r].
\end{equation} 
Even though the existence of a density is not very restrictive in practice, it seems that our results would still be true if we only assumed the existence of a density in the neighborhood $[m-r,m+r]$ of $m$. Moreover, \eqref{comment-ass} is a natural and standard assumption on the distribution of $X_1$ in order to estimate its population median $m$ at the usual $n^{-1/2}$ rate. Indeed, if \eqref{comment-ass} does not hold, then the distribution of $X_1$ does not put enough mass around the median, which becomes harder to estimate. For instance, it is well known that the empirical median of iid random variables is only asymptotically normal when the data have a positive density at the true median (the asymptotic variance being $(4f(m))^{-1}$).

We define the randomized estimator of the median as 
$$\tilde m_T(x)=\hat m_T(x)+\frac{2Z}{\varepsilon}\Sm_h^{(\beta)}(x), x\in\R^n,$$
where $\DS \beta=\frac{\varepsilon}{2\log(2/\delta)}$ and $Z$ is a Laplace random variable with parameter $1$. By Theorem \ref{thm:DP-SS}, $\tilde m_T$ is $(\varepsilon,\delta)$-differentially private, and we have the following theorem. 

\begin{theorem} \label{thm:median-SS}

Let $X=(X_1,\ldots,X_n)$ be a vector of iid random variables satisfying Assumption \ref{Ass:1}. Assume that $|m|\leq R$ for some $R>0$ and let $T>R+r$. Let $\alpha\in [8e^{-nL^2r^2/4},1]$. Then, with probability at least $1-\alpha$,
\begin{align*}
	|\tilde m_T(X)-m| & \leq \sqrt{\frac{2\log\left(\frac{8}{\alpha}\right)}{nL^2}} + \frac{4\log\left(\frac{8}{\alpha}\right)\log\left(\frac{2}{\delta}\right)}{eL\varepsilon^2 n}\left(\log\left\lfloor {Lrn}{2}\right\rfloor+ \log\left(\frac{4}{\alpha}\right)\right) \\ 
	& \hspace{30mm} + \frac{4T\log\left(\frac{4}{\alpha}\right)}{\varepsilon}e^{-\frac{\varepsilon Lrn}{4\log\left(\frac{2}{\delta}\right)}}.
\end{align*}

\end{theorem}

Note that in this theorem, the probability is computed with respect to the joint randomness of the algorithm and of the data. Moreover, our estimator truncates the data but our guarantees do not assume that the data are bounded. Perhaps a drawback of this result is that it assumes that the true median lies in a bounded range.

\begin{proof}
	The proof is decomposed into two parts. First, we bound the smooth sensitivity of $\hat m_T$ evaluated at the random sample $X=(X_1,\ldots,X_n)$ with high probability, using Assumption \ref{Ass:1} (Lemma \ref{lemma:bound-SS} with $\alpha_2=\alpha/4$). Then, again using this assumption, we bound the deviations of the empirical truncated median $\hat m_T(X)$ (Lemma \ref{lemma:bound-emp-med} with $\alpha_3=\alpha/4$). Finally, we get the desired result by using a union bound, where we also control the tails of the Laplace random variable $Z$ and apply the triangle inequality.

\begin{lemma} \label{lemma:bound-SS}

Let $\alpha_2\in (0,1]$. With probability at least $1-\alpha_2-2e^{-nL^2r^2/4}$,
\begin{equation*}
	\Sm_{\hat m_T}^{(\beta)}(X) \leq \frac{1}{eL\beta n}\left(\log\lfloor Lrn/2\rfloor+ \log(1/\alpha_2)\right) +2Te^{-\beta Lrn/2}.
\end{equation*}

\end{lemma}

\begin{lemma} \label{lemma:bound-emp-med}

Let $\alpha_3\in [2e^{-nL^2r^2/2},1]$. Then, with probability at least $1-\alpha_3$,
\begin{equation*}
	|\hat m_T(X)-m|\leq \sqrt{\frac{2\log(2/\alpha_3)}{nL^2}}.
\end{equation*}

\end{lemma}

The restriction on $\alpha_3$ is necessary because Assumption \ref{Ass:1} only imposes a control of the cdf $F$ on a neighborhood of $m$ of size $r$.

\end{proof}

The first term in the upper bound in Theorem \ref{thm:median-SS} is a sub-Gaussian term that comes from the empirical median itself. The other terms, of a smaller order in $n$, are the price to pay in order to apply the Laplace mechanism to the truncated median. Note that $T$ can be chosen as a growing sequence of $n$ so long as the last term does remains very small (while in many applications, $\varepsilon$ may also depend on $n$, making the ratio in the exponential term small).

\subsection{A propose-test-release approach}
\label{sec:MedianPTR}
Our second algorithm is inspired by the propose-test-release paradigm introduced in \cite{dworkandlei2009}. The high level idea of this approach is to propose a bound on the local sensitivity of the desired statistic and test in a differentially private way whether the suggested bound is high enough to ensure privacy. If the proposed bound passes the test then an appropriately calibrated noisy version of the statistic is released. If the proposed bound is not high enough then the algorithm returns ``No Reply", which we denote by $\perp$. Our estimator is carefully calibrated so that the resulting algorithm returns a numerical value with high probability while ensuring satisfactory finite sample statistical guarantees. In particular,  the leading term in the deviations resulting from this procedure coincides with the usual median sub-Gaussian deviations and  dominates the  additional error term  introduced by the privacy inducing mechanism. We require some additional notation in order to define our estimator.  A key component of the algorithm is the quantity

\begin{equation}
\label{Ahat}\hat{A}(x)=\min \left\{k=0,1,\ldots: \exists x'\in\R^n, \dist(x,x')=k, \left|\hat m (x)-\hat{m}(x')\right|>\eta\right\}. 
\end{equation}
It is not too difficult to see that the global sensitivity of $\hat{A}$ equals $1$ since for any $x'$ such that $\dist(x,x')=1$, the only possible values of $\hat{A}(x')$ are $\hat{A}(x)-1,\hat{A}(x)$ and $\hat{A}(x)+1$. Therefore, by Theorem \ref{thm:DP-GS}, $\hat{A}$ can be made $(\varepsilon,0)$-differentially private with the Laplace mechanism
\begin{equation}
\label{Rhat}
\tilde{A}(x)=\hat{A}(x)+\frac{1}{\varepsilon} Z_1,
\end{equation}
where  $Z_1\sim \mathrm{Lap}(1)$.  We introduce another independent random variable $Z_2\sim \mathrm{Lap}(1)$ and define our randomized propose-test-release median estimator as
\begin{equation}
\label{PTR}
\tilde{m}(x)=\begin{cases}
\perp  &\mbox{ if } \tilde A(x) \leq 1+\frac{1}{\varepsilon}\log(2/\delta) \\
\hat{m}(x)+\frac{\eta}{\varepsilon}Z_2 & \mbox{otherwise}
\end{cases},
\end{equation}
for all $x\in\R^n$. Intuitively, $\tilde{m}(x)$ is more likely to output ``No Reply'' for less favorable data configurations $x\in\R^n$ that lead to small values of $\hat{A}(x)$.  Clearly the choice $\eta$ is critical as it controls the size of $\hat{A}(x)$ and hence the probability of no reply, and also the noise term added to $\hat{m}(x)$ when it is released.   

It is interesting to notice that from a computer science perspective $\hat{A}(x)$ can be thought of as the answer to the query: what is the minimum number of observations that we need to change to $x$ before we change the value of the empirical median $\hat{m}(x)$ by at least $\eta$? From a statistical point of view,  it is reminiscent of the finite sample breakdown point studied in robust statistics \cite{donohoandhuber1983,huberandronchetti2009}. While the finite sample breakdown point is usually defined as the minimum number of points that needs to be moved arbitrarily before an estimator becomes infinite, $\hat{A}(x)$ can be interpreted as a relaxed version of the finite sample breakdown point of the median at the scale $\eta$.

\begin{theorem}\label{thm:PTR}
The randomized estimator $\tilde{m}$ is $(2\varepsilon,\delta)$-differentially private and can be computed in $O(n\log(n))$ time. Furthermore, let $X=(X_1,\ldots,X_n)$ be a vector of iid random variables that satisfy Assumption \ref{Ass:1} and let $\alpha\in [8e^{-L^2r^2n/2},1]$ and $\DS C=L^{-1}\left(1+\frac{\log(4/\alpha)}{\log(Lrn/2)}\right)$. Then, choosing $\eta=\frac{C\log(n)}{\varepsilon n}\{\log(2/\delta)+\log(8/\alpha)+\varepsilon\}$, with probability at least $1-\alpha$ we have that 
$$ |\tilde{m}(X)-m| \leq \sqrt{\frac{2\log(8/\alpha)}{nL^2}}+\frac{C\log(n)\{\log(2/\delta)+\log(8/\alpha)+\varepsilon\}\log(8/\alpha)}{\varepsilon^2 n}$$
\end{theorem}
Note that in this theorem (as in Theorem \ref{thm:median-SS}), the probability is computed with respect to the joint randomness of the algorithm and of the data. Furthermore, the second term corresponds to a subexponential type error because the dependence on $\log(1/\alpha)$ does not appear inside a square root.
\begin{proof}
 We only highlight the main ideas of the arguments and relegate some technical details to the proofs of Lemmas \ref{lemma:PTRdp}-\ref{lemma:PTRnoreply}  to the Appendix.  

The estimator $\tilde{m}$ can be shown to be $(2\varepsilon,\delta)$-differentially private by adapting the arguments used for establishing differential privacy of  the propose-test-release median algorithm introduced in \cite{dworkandlei2009} (Lemma \ref{lemma:PTRdp}).  
Its computation is $O(n\log(n))$ time as it comes down to the complexity of sorting $x=(x_1,\dots,x_n)$ and computing $\hat{A}(x)$ using the sorted $x$ for any $x\in\RR^n$.  Once the data is sorted the computation of $\hat{A}(x)$ can be done in $O(n\log(n))$ computations using an interval halving algorithm (Lemma \ref{lemma:PTRcomp}).

For the desired sub-Gaussian deviation inequality we first note that with probability at least $1-\alpha_2-\tau_2$
$$\Big|\hat{m}(X)-m+\frac{\eta}{\varepsilon}Z_2 \Big| \leq \sqrt{\frac{2\log(2/\alpha_2)}{nL^2}}+ \frac{\eta}{\varepsilon}\log(2/\tau_2),$$
where we used the triangle inequality, equation \eqref{boundSS6} in the Appendix and $\PP(|Z_2|>t)=e^{-t}$. It therefore remains to verify that the choice of $\eta$ stated in the theorem suffices to establish that the probability that $\tilde{m}$ gives a ``No Reply'' is at least $\alpha_1+\tau_1$. The main technical obstacle for this is to lower bound $\hat{A}(X)$ with high probability. We establish this result by leveraging  Condition 1 in order to control the maximum gap between consecutive order statistics of $X$ in a neighborhood of $\hat{m}(X)$ (Lemma \ref{lemma:PTRnoreply}).

 \begin{lemma}\label{lemma:PTRdp}
The estimator $\tilde{m}$ is $(2\varepsilon,\delta)-$differentially private.
\end{lemma}

\begin{lemma}\label{lemma:PTRcomp}
 For all $x\in\R^n$, the computation of $\tilde{m}(x)$ is $O(n\log(n))$.
\end{lemma}

\begin{lemma}\label{lemma:PTRnoreply}
Assume that Assumption \ref{Ass:1} holds and let $\eta=\frac{C\log(n)}{\varepsilon n}\{\log(2/\delta)+\log(2/\tau_1)+\varepsilon\}$. Then  for  $\tau_1\in (0,1]$ and $\alpha_1\in(2e^{-\frac{L^2r^2}{2}n},1]$,  we have that $\PP(\tilde{m}(X)=\perp)\leq \tau_1+\alpha_1$.
\end{lemma}

\end{proof}

We note that while the propose-test-release median estimator of \cite{dworkandlei2009} was shown to converge in probability to the population median, it is easy to see that the resulting estimator converges at the rate $n^{1/3}$. Our estimator allows to take smaller values of $\varepsilon$ and $\delta$ while preserving the finite sample sub-Gaussian deviations of the empirical median under minimal conditions.

\section{Discussion on mean estimation}
\label{others}

Recall the definition of the median of means function. Let $K\leq n/2$ be integer and let $N=n/K$, which we assume to be an integer, for the sake of simplicity. We partition the set $\{1,\ldots,n\}$ into $N$ groups $B_1,\ldots,B_N$ of identical size $K$. For all $x=(x_1,\ldots,x_n)\in\R^n$ and each group $B_j, j=1,\ldots,N$, we let $\bar x_j=K^{-1}\sum_{t\in B_j}x_t$, the empirical mean of the coordinates of $x$ indexed in $B_j$. We denote by $\bar x_{(j)}, j=1,\ldots,N$ the empirical means reordered in nondecreasing order. The median of means of $x$ is defined as $\hat\mu(x)=\bar x_{(\ell)}$, where now, $\ell=\lfloor N/2\rfloor$. 

It is well known that if $X_1,\ldots,X_n$ are iid random variables with only the first two moments assumed to be finite $\mu=\E[X_1]$ and $\sigma^2=\Var[X_1]$, then the median of mean estimator $\hat\mu(X_1,\ldots,X_n)$ satisfies the following, for $N=8\log(2/\alpha)$ and $K=n/N$, where $\alpha\in (0,1)$ is a prescribed probability level.

\begin{lemma}\cite[Lemma 2]{bubecketal2013}
	With probability at least $1-\alpha$, 
$$|\hat\mu(X_1,\ldots,X_n)-\mu|\leq 2\sigma\sqrt{\frac{\log(2/\alpha)}{n}}.$$
\end{lemma}

However, the Laplace mechanism that we have used in Section \ref{sec:MedianSS} cannot be applied directly to $\hat\mu$. It is easy to see that the global sensitivity $\GS_{\hat \mu}$ is infinite. The local sensitivity is finite as long as the empirical means $\bar x_{(j)}, j=1,\ldots,N$ take at least three different values. In that case, $\LS_{\mu}(x)=\max\left(\bar x_{(\ell)}-\bar x_{(\ell-1)},\bar x_{(\ell+1)}-\bar x_{(\ell)}\right)$. However, the $\beta$-smooth sensitivity $\Sm_{\hat\mu}^{(\beta)}(x)$ is infinite, for all $x\in\R^n$.

A natural alternative to $\hat\mu$ consists of truncating the empirical means $\bar x_j, j=1,\ldots,N$, before taking their median, or to first truncate the $x_i$'s, compute the new empirical means on each block and take the median. These two functions are different, but we only focus on the first one, since a similar analysis would hold for the second one (with a different threshold).

For $T>0$, we define $\hat\mu_T(x)=f_T(\bar x_{(\ell)})$.

\begin{lemma} \label{lemma:SS-MedMeans}
	Let $\beta>0$. Then, the $\beta$-smooth sensitivity of the function $\hat \mu$ is given by
\begin{equation*}
	\Sm_{\hat \mu_T}^{(\beta)}=\max_{0\leq j\leq N}e^{-\beta j} \max_{t=0,\ldots,j+1} \left(f_T(\bar x_{(\ell+t)})-f_T(\bar x_{(\ell+t-j-1)})\right),
\end{equation*}
where we set $\bar x_{(j)}=\infty$ if $j\geq N+1$ and $\bar x_{(j)}=-\infty$ if $j\leq -1$.
Moreover, it is bounded from above by $2T$.
\end{lemma}

The proof of the first this lemma is similar to that of Lemma \ref{lemma:bound-SS} and is omitted here. For the second part, it is clear that one can bound $\Sm_{\hat \mu_T}^{(\beta)}$ by $f_T(\bar x_{(N)})-f_T(\bar x_{(1)})$ and this bound is tight up to a constant. In general, unless $T$ is significantly larger than $\sqrt n$, $\bar X_{(N)}\geq T$ and $\bar X_{(1)}\leq -T$ with high probability. Therefore, the Laplace mechanism that we used in \ref{sec:MedianSS} will not provide reasonable deviation bound. Actually, the same limitations also hold for the truncated mean estimator studied in \cite{caietal2019}, where the authors deal with a truncation parameter $T$ of order $\log n$ since their data are sub-Gaussian, whereas under the current assumptions, $T$ should be taken of the order of $\sqrt n$. In general, it seems that the usual estimators of the mean under the existence of only low moments (median of means, truncated mean, Catoni's estimators \cite{catoni2012,bubecketal2013,devroyeetal2016}) can not be adapted in a straightforward way in order to obtain differentially private estimators that admit a sub-Gaussian error plus a term of smaller order.

The propose-test-release approach  introduced in \ref{sec:MedianPTR} will also fail to give sub-Gaussian type deviations if applied to the median of means function. Indeed in order to obtain inequalities similar to the ones obtained in Theorem \ref{thm:PTR} one would have to choose the constant $\eta$ in equation \eqref{Ahat} to be of the order $1/n$. However one can show that the expect the empirical means $\bar X_1,\dots,\bar X_n$ are separated by a distance of at least of the order $1/\sqrt{n\alpha}$ with probability $1-\alpha$. Consequently if we choose $\eta$ smaller than $1/\sqrt{n\alpha}$ the probability of ``No Reply'' of the propose-test-release approach will be too large since we will always have small values of $\hat A(x)$. A similar argument suffices to see that our propose-test-release mechanism will also fail to give the desired deviations for the truncated mean estimator studied in \cite{caietal2019}.

\section{Conclusion}

We studied the problem of differentially private location parameter estimation from a non-asymptotic deviations perspective.  In particular, we proposed the first private median estimators that exhibit a leading sub-Gaussian error  terms with high probability. Our first estimator uses truncation in order to get high probability control of the smooth sensitivity while the second one avoids truncating the data by first identifying in a privacy-preserving fashion whether the data are in a favorable configuration before releasing a noisy version of the desired output.

We showed that differentially private versions of well known sub-Gaussian mean estimators fail to exhibit the optimal deviations of their nonprivate counterparts. A possible explanation for these weaker results for mean estimation is that perhaps one can only expect to obtain sub-Gaussian type deviations for  differentially private versions of robust statistics in the sense of \cite{huberandronchetti2009}. It would therefore be interesting to extend our methods to other robust statistics beyond the median.

\bibliographystyle{plain}
\bibliography{Biblio}

%
%
%
%
%
%

\appendix
\section*{Appendix: Differentially private sub-Gaussian location estimators}

\subsection*{Proof of Lemma \ref{lemma:bound-SS}}

Let $k_0=\left\lfloor\frac{Lrn}{2}\right\rfloor$. Note that $k_0\leq n/4$, since $2Lr\leq F(m+r)-F(m-r)\leq  1$, by Assumption \ref{Ass:1}. We decompose $\Sm_h^{(\beta)}(X)$ into two parts:
\begin{align}
	\Sm_{\hat m_T}^{(\beta)}(X) & = \max_{k=0,\ldots,k_0-1}e^{-\beta k}\max_{t=0,\ldots,k+1}\left(Y_{(\ell+t)}-Y_{(\ell+t-k-1)}\right) \nonumber \\
	& \quad \quad + \max_{k\geq k_0}e^{-\beta k}\max_{t=0,\ldots,k+1}\left(Y_{(\ell+t)}-Y_{(\ell+t-k-1)}\right) \label{boundSS1}
\end{align}
The second term is upper bounded by $2Te^{-\beta k_0}$. For the first term, we consider the event $\DS \mathcal A=\{Y_{(i+1)}-Y_{(i)}\leq \frac{C\log k_0}{n}, \forall i=\ell-k_0,\ldots,\ell+k_0-1\}$, where $\DS C=L^{-1}\left(1+\frac{\log(1/\alpha_2)}{\log(k_0)}\right)$. If $\mathcal A$ is satisfied, the first term in \eqref{boundSS1} is bounded by 
\begin{align*}
	\max_{0\leq k\leq k_0-1} e^{-\beta k}k\frac{C\log k_0}{n} & = \max_{0\leq k\leq k_0-1} e^{-\beta k}\beta k\frac{C\log k_0}{\beta n} \\
	& \leq \left(\max_{u\geq 0}e^{-u}u\right)\frac{C\log k_0}{\beta n} = \frac{C\log k_0}{e\beta n},
\end{align*} 
hence, if $\mathcal A$ is satisfied, \eqref{boundSS1} yields that 
\begin{equation} \label{boundSS2}
	\Sm_{\hat m_T}^{(\beta)}(X) \leq \frac{C\log k_0}{e\beta n}+2Te^{-\beta k_0}.
\end{equation}
Now, we control the probability of the complement of $\mathcal A$. Consider the event $\mathcal B=\{m-r<X_{(\ell-k_0)}\leq X_{(\ell+k_0)}<m+r\}$.
By Assumption \ref{Ass:1}, if $\mathcal B$ is satisfied, it holds that for all $i=\ell-k_0,\ldots,\ell+k_0$, $X_{(i)}=Y_{(i)}$, and for all $i=\ell-k_0,\ldots,\ell+k_0-1$,
\begin{align*}
	F(X_{(i+1)})-F(X_{(i)}) & \geq L(X_{(i+1)}-X_{(i)}) \\
	& = L(Y_{(i+1)}-Y_{(i)}),
\end{align*}
hence, if $\mathcal A^{\complement}$ holds together with $\mathcal B$, then 
\begin{align*}
	\max_{i=\ell-k_0,\ldots,\ell+k_0-1} U_{(i+1)}-U_{(i)}\geq \frac{CL\log k_0}{n},
\end{align*}
where we have denoted by $U_j=F(X_j), j=1,\ldots,n$. Therefore, 
\begin{align}
\PP[\mathcal A^{\complement}] & = \PP[\mathcal A^{\complement}\cap \mathcal B] + \PP[\mathcal A^{\complement}\cap \mathcal B^{\complement}] \nonumber \\
& \leq \PP\left[\max_{i=\ell-k_0,\ldots,\ell+k_0-1} U_{(i+1)}-U_{(i)}\geq \frac{CL\log k_0}{n}\right] +\PP[\mathcal B^{\complement}] \nonumber \\
& \leq \sum_{i=\ell-k_0}^{\ell+k_0-1}\PP\left[U_{(i+1)}-U_{(i)}\geq \frac{CL\log k_0}{n}\right]+\PP[\mathcal B^{\complement}], \label{boundSS3}
\end{align}
where the last inequality follows from a union bound. By Assumption \ref{Ass:1}, the $U_j$'s are iid and uniformly distributed in $[0,1]$. Therefore, for all $i=1,\ldots,n-1$, $U_{(i+1)}-U_{(i)}$ has the same distribution as $U_{(1)}$ and Therefore, for all $i=\ell-k_0,\ldots,\ell+k_0-1$,
\begin{align}
	\PP\left[U_{(i+1)}-U_{(i)}\geq \frac{CL\log k_0}{n}\right] & = \PP\left[U_{(1)}\geq \frac{CL\log k_0}{n}\right] \nonumber \\
	& = \PP\left[U_j\geq \frac{CL\log k_0}{n}, \forall j=1,\ldots,n\right] \nonumber \\
	& = \PP\left[U_1\geq \frac{CL\log k_0}{n},\right]^n \nonumber \\
	& = \left(1-\frac{CL\log k_0}{n}\right)^n \leq e^{-CL\log k_0}. \label{boundSS4}
\end{align}
In order to bound the probability of $\mathcal B^{\complement}$, it suffices to note that 
\begin{equation} \label{boundSS5}
	\PP[\mathcal B^{\complement}] = \PP[X_{(\ell+k_0)}>m+r]+\PP[X_{(\ell-k_0)}<m-r].
\end{equation}
We only bound the first term in the right hand side of \eqref{boundSS5}, since the second term can be bounded with similar computations. By definition of the order statistics,
\begin{align}
	& \PP[X_{(\ell+k_0)}>m+r] \nonumber \\
	& \quad = \PP \left[\sum_{i=1}^n \mathds 1_{X_i>m+r} \geq n-(\ell+k_0)+1\right] 
	\nonumber \\
	& \quad = \PP\left[\sum_{i=1}^n \left(\mathds 1_{X_i>m+r}-(1-F(m+r))\right) \geq n-(\ell+k_0)+1-n(1-F(m+r))\right]
	 \nonumber \\
	& \quad \leq \exp\left(-\frac{n}{2}\left(F(m+r)-\frac{\ell+k_0}{n}\right)^2\right), \label{boundSS6}
\end{align}
where the last inequality follows from Hoeffding's inequality. By Assumption \ref{Ass:1}, $\DS F(m+r)-\frac{\ell+k_0}{n}\geq F(m)+Lr-\frac{1}{2}-\frac{k_0}{n} = Lr-\frac{k_0}{n}\geq Lr/2>0$. Thus, using \eqref{boundSS6}, we have that $\DS \PP[X_{(\ell+k_0)}>m+r] \leq e^{-nL^2r^2/4}=: \alpha_1$. Finally, we get that $\DS \PP[\mathcal B^{\complement}]\leq 2\alpha_1$. So, by \eqref{boundSS3}, 
\begin{align*}
	\PP[\mathcal A^{\complement}] & \leq 2k_0e^{-CL\log k_0}+2\alpha_1 \\
	& = 2k_0^{1-CL}+2\alpha_1 = \alpha_2+2\alpha_1.	
\end{align*}

This ends the proof of Lemma \ref{lemma:bound-SS}.

\subsection*{Proof of Lemma \ref{lemma:bound-emp-med}}

Let $t\in [0,r]$. Then, 
\begin{equation} \label{boundemp1}
	\PP[|\hat m(X)-m|>t] = \PP[\hat m(X)>m+t]+\PP[\hat m(X)<m-t].
\end{equation}
Now, we only bound the first term, since the second term would be treated with the exact same arguments. By definition of the empirical median,
\begin{align*}
	\PP[\hat m(X)>m+t] & = \PP[\sum_{i=1}^n \mathds 1_{X_i>m+t}\geq n/2] \\
	& = \PP[\sum_{i=1}^n \left(\mathds 1_{X_i>m+t}-(1-F(m+t))\right)\geq n\left(F(m+t)-1/2\right)] \\
	& \leq \exp\left(-\frac{n}{2}(F(m+t)-1/2)^2\right) \\
	& \leq e^{-nL^2t^2/2},
\end{align*}
where we used Hoeffding's inequality and Assumption \ref{Ass:1}. Therefore, we simply get, from \eqref{boundemp1}, that 
\begin{equation} \label{boundemp3}
	\PP[|\hat m(X)-m|>t] \leq 2e^{-nL^2t^2/2}.
\end{equation}
Now, since $-T<-R-r\leq m-r\leq m-t\leq m+t\leq m+r\leq R+r<T$, if $|\hat m(X)-m|\leq t$, then it must hold that $\hat m(X)=\hat m_T(X)$, therefore, \eqref{boundemp3} yields that 
\begin{equation*} 
	\PP[|\hat m_T(X)-m|>t] \leq \PP[|\hat m(X)-m|>t] \leq 2e^{-nL^2t^2/2}.
\end{equation*}
This concludes the proof of the lemma, by taking $\DS t=\sqrt{\frac{2\log(2/\alpha_3)}{nL^2}}$.

\subsection*{Proof of Lemma \ref{lemma:PTRdp}}
Since  $\tilde{m}(x)=\perp \iff \tilde{A}(x)\leq 1+\frac{1}{\varepsilon}\log(2/\delta)$ and  $\tilde{A}$ is $(\varepsilon,0)$-differentially private, it follows that 
\begin{equation} \label{lemma4:noreply} 
	\PP[\tilde{m}(x)=\perp]\leq e^\varepsilon \PP[\tilde{m}(x')=\perp],
\end{equation} 
for all $x,x'$ such that $\dist(x,x')\leq 1$. Indeed outputting ``No Reply'' has the same privacy guarantee as $\tilde A$ because differential privacy is not affected by post-processing \cite[Proposition 2.1]{dworkandroth2014}.

Let $x,x'\in\R^n$ with $\dist(x,x')\leq 1$. We will now show that for all Borel sets $\mathcal O\subseteq \R$,
\begin{equation}
\label{lemma4:reply}
\PP[\tilde{m}(x)\in\mathcal{O}]\leq  e^{2\varepsilon}\PP[\tilde{m}(x')\in \mathcal{O}]+\delta.
\end{equation}
Note that for $\tilde m(x)$ to be a real number, it has to be that the estimator has outputted a reply, i.e., $\tilde A(x)>1+(1/\varepsilon)\log(2/\delta)$. On the one hand, if  $|\hat{m}(x)-\hat{m}(x')|\leq \eta$, we have 
\begin{align}
\label{D0}
 \PP[\tilde{m}(x)\in\mathcal{O}]&=\PP\Big[\hat{m}(x)+\frac{\eta}{\varepsilon}Z_2\in\mathcal{O},~\hat{A}(x)+\frac{1}{\varepsilon}Z_1> 1+\frac{1}{\varepsilon}\log(2/\delta)\Big] \nonumber \\
 & =\PP\Big[\hat{m}(x)+\frac{\eta}{\varepsilon}Z_2\in\mathcal{O}\Big]\PP\Big[\hat{A}(x)+\frac{1}{\varepsilon}Z_1 >1+\frac{1}{\varepsilon}\log(2/\delta)\Big]\nonumber \\
  &\leq e^\varepsilon \PP\Big[\hat{m}(x')+\frac{\eta}{\varepsilon}Z_2\in\mathcal{O}\Big]e^\varepsilon\PP\Big[\hat{A}(x')+\frac{1}{\varepsilon}Z_1 >1+\frac{1}{\varepsilon}\log(2/\delta)\Big] \nonumber\\
 &=e^{2\varepsilon}\PP[\tilde{m}(x')\in \mathcal{O}]\nonumber \\
 &\leq e^{2\varepsilon}\PP[\tilde{m}(x')\in \mathcal{O}]+\delta,
\end{align}
by the sliding property of the Laplace distribution \cite[Section 2.1.1]{nissimetal2007} and where the second and the last equalities used independence of $Z_1$ and $Z_2$.  On the other hand, if $|\hat{m}(x)-\hat{m}(x')| > \eta$ then $\hat{A}(x)=\hat{A}(x')=1$ which in turn entails that
\begin{align}
\label{notD0}
\PP[\tilde{m}(x)\in\mathcal{O}]&=\PP\Big[\hat{m}(x)+\frac{\eta}{\varepsilon}Z_2\in\mathcal{O},~\frac{1}{\varepsilon}Z_1> \frac{1}{\varepsilon}\log(2/\delta)\Big]\nonumber \\
&\leq \PP\Big[\frac{1}{\varepsilon}Z_1 >\frac{1}{\varepsilon}\log(2/\delta)\Big]\nonumber\\
&=\delta \nonumber \\
&\leq e^{2\varepsilon}\PP[\tilde{m}(x')\in\mathcal{O}]+\delta
\end{align}
Therefore, combining \eqref{D0} and \eqref{notD0} yields \eqref{lemma4:reply}. Now, let $\mathcal O'$ be a Borel set of the extended real line $\R\cup\{\perp\}$. Then, $\mathcal O'$ is equal to either $\mathcal O$, or $\mathcal O\cup\{\perp\}$, for some Borel set of $\R$. In the former case, \eqref{lemma4:reply} concludes the proof of the lemma. In the latter case, we write, for all $x,x'\in\R^n$ with $\dist(x,x')\leq 1$:
\begin{align*}
	\PP[\tilde{m}(x)\in\mathcal{O'}] & = \PP[\tilde{m}(x)\in\mathcal{O}] + \PP[\tilde{m}(x)=\perp] \\
	& \leq e^{2\varepsilon}\PP[\tilde{m}(x')\in\mathcal{O}]+\delta + e^\varepsilon \PP[\tilde{m}(x')=\perp] \\
	& \leq e^{2\varepsilon} \PP[\tilde{m}(x')\in\mathcal{O'}]+\delta,
\end{align*}
thanks to \eqref{lemma4:noreply} and \eqref{lemma4:reply}.

\subsection*{Proof of Lemma \ref{lemma:PTRcomp}}

It suffices to show that $\hat{A}(x)$ can be computed in near linear time.  For this we can  sort $x$ and take the resulting order statistics $x_{(1)},\dots,x_{(n)}$ to  compute $\hat{A}(x)=\min\{k \mbox{ s.t. } \max_{0\leq t\leq k+1} ( x_{(m+t)}-x_{(m+t-k-1)})>\eta\}$.  For a fixed $k$, solving $\max_{0\leq t\leq k+1} ( x_{(m+t)}-x_{(m+t-k-1)})$ takes at most $O(n)$ operations. Furthermore, using a dichotomy method (which is valid since $\max_{0\leq t\leq k+1} ( x_{(m+t)}-x_{(m+t-k-1)})$ is monotone in $k$), we see that we only need to explore $O(\log(n))$ values of $k$ in order to find $\hat{A}(x)$. Hence for sorted $x$, we showed that $\hat{A}(x)$  can be computed in in $O(n\log (n))$ time. Since the initial sorting step also takes $O(n \log (n))$ operations, the overall algorithm is $O(n\log(n))$.

\subsection*{Proof of Lemma \ref{lemma:PTRnoreply}}
 Using the notation that we used in the proof of Lemma \ref{lemma:bound-SS}, it follows from the definition of $\tilde{m}(X)$ that
\begin{align*}
\PP[\tilde{m}(X)=\perp] &=\PP\bigg[\hat{A}(X)+\frac{1}{\varepsilon}Z_1 \leq 1+\frac{1}{\varepsilon}\log(2/\delta)\bigg] \\
&\leq \PP\bigg[\Big\{\hat{A}(X)+\frac{1}{\varepsilon}Z_1 \leq 1+\frac{1}{\varepsilon}\log(2/\delta)\Big\}\cap \mathcal{B}\bigg] +\PP[\mathcal B^{\complement}]\\
 & \leq \PP\bigg[ \Big\{\frac{\eta n}{C\log(n)}+ \frac{1}{\varepsilon}Z_1\leq 1+\frac{1}{\varepsilon}\log(2/\delta)\Big\}\cap \mathcal B\bigg]  +\alpha_1\\
 &  \leq \PP\bigg[\frac{\eta n}{C\log(n)} + \frac{1}{\varepsilon}Z_1\leq1+\frac{1}{\varepsilon}\log(2/\delta)\bigg] +\alpha_1 \\
 &=\PP\bigg[Z_1\leq \varepsilon\Big(1+ \frac{1}{\varepsilon}\log(2/\delta)-\frac{\eta n}{C\log(n)}\Big)\bigg]+ \alpha_1.
 \end{align*}
In order to set the above bound equal to $\tau_1+\alpha_1$ we have to use that $\PP[Z_1\leq -\log(2/\tau_1)]=\tau_1$. This requires that  
$$ -\frac{1}{\varepsilon}\log(2/\tau_1)=1+ \frac{1}{\varepsilon}\log(2/\delta)-\frac{\eta n}{C\log(n)}$$
from which we obtain that
$$ \eta=\frac{C\log(n)}{\varepsilon n}\Big(\log(2/\delta)+\log(2/\tau_1)+\varepsilon\Big). $$

\end{document}